\documentclass[12pt]{amsart}
\usepackage{amssymb}
\usepackage{graphics}
\usepackage{amsmath}
\usepackage{enumerate}
\newcommand{\Z}{{\mathbb Z}}

\newcommand{\Q}{{\mathbb Q}}

\theoremstyle{plain}

\newtheorem{theorem}{Theorem}

\newtheorem{corollary}{Corollary}

\theoremstyle{definition}
\newtheorem{definition}{Definition}

\begin{document}
\title[Mild pro-$p$-groups with $4$ generators]
     { Mild pro-$p$-groups with $4$ generators}
\author{Michael R. Bush}
\address{Dept. of Mathematics \& Statistics,
  University of Massachusetts, Amherst, MA 01003-9305, USA}
\email{bush@math.umass.edu}
\author{John Labute}
\address{Dept. of Mathematics \& Statistics,
  McGill University, Burnside Hall, 805 Sherbrooke Street West,
  Montreal, QC H3A 2K6, Canada}
\email{labute@math.mcgill.ca}
\date{December 20, 2005}
\begin{abstract}
Let $p$ be an odd prime and $S$ a finite set of primes $\equiv1$
mod $p$. We give an effective criterion for determining when the
Galois group $G=G_S(p)$ of the maximal $p$-extension of $\mathbb
Q$ unramified outside of $S$ is mild when $|S|=4$ and the cup
product $H^1(G,\Z/p\Z)\otimes H^1(G,\Z/p\Z) \rightarrow
H^2(G,\Z/p\Z)$ is surjective.
\end{abstract}

\maketitle

\section{Introduction}
Let $p$ be an odd prime and $S$ a finite set of primes not
containing $p$. Let $G=G_S(p)$ be the Galois group of the maximal
$p$-extension of $\mathbb{Q}$ unramified outside $S$. We can
assume $S = \{ q_1,\ldots\,,q_m \}$ with $q_i \equiv 1 \mod p$.
Work of Koch \cite{Ko}  shows that $G = F/R$ where $F$ is the free
pro-$p$-group on $x_1,\ldots\,,x_m$ and $R = (r_1,\ldots\,,r_m)$
with
\[  r_i \equiv x_i^{q_i - 1} \prod_{j \neq i} [x_i,x_j]^{l_{i j}} \mod F_3 \]
where $l_{i j} \in \mathbb{F}_p$ and $F_n$ is the n-th term in the
lower $p$-central series defined recursively by $F_1 = F$ and
$F_{n+1} = F_n^p [F_n,F]$. Moreover, $l_{i j}$ is the image in
$\mathbb{F}_p$ of any integer $r$ satisfying
$$
q_i\equiv g_j^{-r} \text{ mod }q_j
$$
where $g_j$ is a primitive root for the prime $q_j$. If
$\chi_1,\ldots,\chi_m\in H^1(G,\Z/p\Z)$ with
$\chi_i(x_j)=\delta_{ij}$, we have
$\chi_i\cup\chi_j(r_i)=\ell_{ij}$, after identifying
$H^2(G,\Z/p\Z)$ with the dual of $R/R^p[R,F]$ via the
transgression map (see \cite{NSW}, Proposition 3.9.13). It follows
that the cup product
$$
H^1(G,\Z/p\Z)\otimes H^1(G,\Z/p\Z) \rightarrow
H^2(G,\Z/p\Z)
$$
is surjective if and only if the images of $r_1,\ldots,r_m$ in
$F_2/F_1^pF_3$ are linearly independent. The latter is true for
any minimal presentation $\langle x_1,\ldots\,,x_m \mid
r_1,\ldots\,,r_d \rangle$ of a pro-$p$-group $G$. The presentation
is said to be of {\em Koch type} if $d\le m$ and the relations
$r_i$ satisfy a congruence of the form
\[  r_i \equiv x_i^{p a_i} \prod_{j \neq i} [x_i,x_j]^{a_{i j}} \mod F_3. \]

In \cite{La} the second author has shown that under certain
conditions on the relation set $R$ of a presentation  the
associated pro-$p$-group $G = F/R$ has many nice properties. These
conditions can often be shown to hold when the presentation is of
Koch type even if the exact form of the relations is undetermined.
This is of particular interest in the case where $G = G_S(p)$. We
recall the main definitions and some of the results here for the
reader's convenience.

Let $G$ be a pro-$p$-group. The lower $p$-central series $G_n$
(defined above) can be used to construct a graded
$\mathbb{F}_p$-vector space $\text{gr}(G) = \bigoplus_{n \geq 1} \text{gr}_n(G)$ where
$\text{gr}_n(G) = G_n / G_{n+1}$. This has the additional
structure of a Lie algebra over the polynomial ring
$\mathbb{F}_p[\pi]$ where multiplication by $\pi$ is induced by
the map $x \mapsto x^p$ and the bracket operation by the
commutator operation in $G$.

Now suppose that $G = F/R = \langle x_1,\ldots\,,x_m \mid
r_1,\ldots\,,r_d
\rangle$ is finitely presented. If $\xi_i$ is the image of $x_i$ in
$\text{gr}_1(F)$ then $\text{gr}(F)$ is the free Lie algebra on
$\xi_1,\ldots\,,\xi_m$ over $\mathbb{F}_p[\pi]$. We let $h_i$
denote the largest value of $n$ for which $r_i \in F_n$ and
let $\rho_i \in \text{gr}_{h_i}(F)$ be the image of $r_i$ under the
canonical epimorphism. We call $\rho_i$ the {\em initial form} of
$r_i$. If $\mathfrak{r}$ is the ideal of $L = \text{gr}(F)$
generated by $\rho_1,\ldots\,,\rho_d$ and $\mathfrak{g} =
L/\mathfrak{r}$ then $\mathfrak{r}/[\mathfrak{r},\mathfrak{r}]$ is
a module over the enveloping algebra $U_\mathfrak{g}$ of
$\mathfrak{g}$ via the adjoint representation.

\begin{definition}
The sequence $\rho_1,\ldots\,,\rho_d$ with $d \geq 1$ is said to be {\em
strongly free} if $U_\mathfrak{g}$ is a free $\mathbb{F}_p[\pi]$-module and
$M = \mathfrak{r}/[\mathfrak{r},\mathfrak{r}]$ is a free $U_\mathfrak{g}$-module
on the images of $\rho_1,\ldots\,,\rho_d$ in $M$. If a pro-$p$-group $G$ has
a finite presentation $F/R$ in which the initial forms of the relators form a strongly
free sequence then $G$ will be called {\em mild}.
\end{definition}

Mild groups $G$ enjoy many nice properties. In particular, the
graded Lie algebra $\text{gr}(G)$ is finitely presented with
presentation $L/\mathfrak{r}$,  the Poincare series for the
enveloping algebra of $\text{gr}(G)$ is given by the formula
\[ P(t) = \frac{1}{(1-t)(1 - mt + t^{h_1} + \ldots + t^{h_d})}\]
and $G$ has cohomological dimension $2$, (cf. \cite{La}, Theorem
2.1).

In \cite{La}, Theorem 3.3, a criterion for determining strong
freeness is given. It amounts to a certain independence condition
on the sequence $\rho_1,\ldots\,,\rho_d$. With this condition one
can easily generate presentations which yield mild pro-$p$-groups.
One example of particular importance is the cycle presentation
with $n$ generators $x_1,\ldots\,,x_n$ and $n$ relators
$[x_1,x_2], [x_2,x_3],\ldots\,,[x_n,x_1] \in \text{gr}_2(F)$. The criterion
also makes it possible to show that the Galois groups $G_S(p)$ are
mild for various choices of $S$ and $p$.

One issue that arises immediately is the variability in the
applicability of the criterion among presentations for the same
group $G$. It is possible for a group that is mild to have
presentations which cannot be shown to be strongly free using the
criterion mentioned above. How then does one recognize whether or
not a group is mild given only one particular presentation?

In this paper we consider this question in the case where $m =
d=4$ and the initial forms of the relators are quadratic (see below). By
\cite{La}, Theorem 3.10, a sequence of initial forms
$\rho_1,\ldots\,,\rho_d$ is strongly free if and only if
$\overline{\rho_1},\ldots\,,\overline{\rho_d}$ is strongly free
where $\overline{\rho_i}$ is the image of $\rho_i$ in the free
$\mathbb{F}_p$-Lie algebra $\overline{L} = L/\pi L$. Moreover, by
\cite{La}, Theorem 3.2, if $\overline{\mathfrak r}$ is the ideal
of $\overline{L}$ then the sequence
$\overline{\rho_1},\ldots\,,\overline{\rho_d}$ is strongly free if
and only if the Poincar\'e polynomial of the enveloping algebra of
$\overline{L}/\overline{\mathfrak r}$ is
$$
\overline{P}(t)=\frac{1}{1-mt+mt^2}.
$$
Note also that $\overline{\rho_1},\ldots\,,\overline{\rho_d}$ are
linearly independent if and only if the cup product
$$
H^1(G,\Z/p\Z)\otimes H^1(G,\Z/p\Z) \rightarrow H^2(G,\Z/p\Z)
$$
is surjective.

By the remarks above we may  work over the field $\mathbb{F}_p$. In fact the
results we shall prove hold more generally so we let $k =
\mathbb{F}_q$ where $q = p^r$ with $r > 0$. From now on $L$ will
 denote the free Lie algebra on $X = \{x_1,\ldots\,,x_4\}$ over $k$.
The Lie algebra $L$ has the usual grading $\bigoplus_{n = 1}^\infty L_n$ obtained by
assigning a weight of $1$ to each generator in $X$. In this grading $L_1$
is the $4$-dimensional $k$-vector space with basis $X$. The component $L_2$ is 6-dimensional with
basis the images of the brackets of pairs of generators of $L$.  The element
$[x_i,x_j]$ will be denoted $x_{ij}$ and by abuse of
notation this will also be used to denote its image in any quotient of $L$.

We will be interested in finitely presented Lie algebras
$L/\mathfrak{r}$ where the ideal $\mathfrak{r}$ is generated by a
set $R$ of 4 relators  $\rho_1,\ldots\,\rho_4 \in L_2$ which are
linearly independent over $k$. We will call such a Lie algebra
{\em quadratic of relation rank 4}. The quadratic algebra
$L/\mathfrak{r}$ is said to be of Koch type if
$\rho_i=\sum_jl_{ij}x_{ij}$. Any automorphism of $L$ that respects
the grading maps a set of quadratic relations $R$ into another
such set $R'$. It is clear that the sequence $R$ is strongly free
if and only if this is true for $R'$. Our main result will be to
show that under the action of a particular group of
transformations there are exactly $4$ equivalence classes of such
sets of relations, two of which are strongly free and two of which
are not.

The proof is constructive and yields a procedure for actually recognizing
which of the $4$ classes contains any given presentation and in particular
whether or not it is mild. It turns out that two such quadratic algebras are
isomorphic if and only if they are isomorphic modulo the $5$-th term of
their lower central series.

\section{Orbits of Presentations}
We fix the lexicographic ordering  $12 < 13 < 14 < 23 < 24 < 34$
on the set of basis elements $\{x_{ij}\}_{1\le  i<j\leq 4}$ of $L_2$.
Any ordered set $R$ of $4$ quadratic relations is now represented by
a $4 \times 6$ matrix in the obvious way. We have two natural
group actions on the space of $4 \times 6$ matrices over $k$.
A left action by $GL_4(k)$ defined by left multiplication and a right
action by $GL_4(k)$ defined by right multiplication after applying the
homomorphism  $\psi : GL_4(k) \cong \text{Aut}(L_1)  \rightarrow
\text{Aut}(L) \rightarrow \text{Aut}(L_2) \cong GL_6(k)$ (the first map is
a lift using the freeness of $L$, and the second map is restriction).
More explicitly if $A \in GL_4(k) \cong \text{Aut}(L_1)$ is defined by
$ A x_i = \sum_{j = 1}^4 a_{j i} x_j $ with $a_{j i} \in k$,
then $\hat{A} = \psi(A) \in GL_6(k)$ satisfies
\[ \hat{A} x_{i j} = \sum_{r s} (a_{r i} a_{s j} - a_{s i} a_{r j}) x_{r s} \]
where the summation is over all pairs $r s$ ordered lexicographically as
described above. Thus $\psi$ viewed as a representation of $GL_4(k)$
is simply the exterior square $\bigwedge^2(k^4)$.

The left and right  actions are compatible and give rise to various transformations
of the corresponding presentations.
It is clear that the isomorphism type of the Lie algebra associated to a
presentation is preserved under both of these actions.
We would like to understand the orbit decomposition
under this double action. The orbits under the left action
correspond to $4$-dimensional subspaces of $k^6$.
Thus we are reduced to understanding the right action
by the group $G = \psi(GL_4(k))$
on the space of all 4-dimensional subspaces of $k^6$. The problem
could therefore be formulated as determining the orbits under the
action of $G$ (or its image in $PGL_6(k)$) on the Grassmanian
$Gr_k(6,4)$.

Note that throughout this paper we will represent subspaces of $k^6$
by listing a basis of row vectors usually in the form of a matrix.

A generating set for the group $G$ can be obtained by applying
$\psi$ to a generating set for $GL_4(k)$. The images of all
elementary matrices form such a set. We introduce some notation to
describe these. Fix an $n \times n$ identity matrix. Let $E^a_{i
j}$ be the elementary matrix obtained by applying the column
transformation $c_j \rightarrow c_j + a c_i$ for $a \in k$. Let
$E_{i j}$ be the elementary matrix obtained by swapping columns
$c_i \leftrightarrow c_j$. Let $E_i^a$ be the elementary matrix
obtained by re-scaling the $i$-th column $c_i \rightarrow a c_i$.
We now introduce notation for the images of such matrices under
$\psi$. For $i \neq j \in \{1,\ldots,\,4\}$ and $a \in k^\times$
let $M_{i j}^a = \psi (E_{i j}^a$) , let $T_{i j} = \psi(E_{i j})$
and let $S_i^a = \psi(E_i^a)$ (where $a \in k^\times$).

There is one additional simplification which we wish to make.
Every 4-dimensional subspace $U$ has a unique orthogonal complement
$U^\perp$ of dimension 2 with respect to the standard inner product on $k^6$.
There is an induced right action of $G$ on the space of 2-dimensional subspaces given
by $U^\perp \cdot M = U^\perp (M^{-1})^{Tr}$. One can show easily that the group $G$ is
closed under taking transposes so our problem is equivalent to understanding the
orbits of 2-dimensional subspaces under right multiplication by elements of $G$.

We are now ready to start investigating the action of $G$. As an intermediate step
we investigate its action on 1-dimensional subspaces. We have the following result.

\begin{theorem}
The space of 1-dimensional subspaces of $k^6$ under the action of
$G$ splits into at most two orbits. Equivalently every 1-dimensional
space is equivalent under $G$ to either
$\begin{bmatrix} 0 & 0 & 0 & 0 & 0 & 1\end{bmatrix}$ or
$\begin{bmatrix} 0 & 0 & 1 & 1 & 0 & 0\end{bmatrix}$.
\end{theorem}
\begin{proof}
Start with an arbitrary nonzero vector
$\begin{bmatrix} * & * & * & * & * & *\end{bmatrix}$.
In general we will use $*$ to represent
an arbitrary field element in a vector to avoid introducing large numbers  of indeterminates.
To refer to a particular component we will use a subscript $*_i$. Let us focus on the
first three entries, there are 3 cases.
\begin{itemize}
\item $*_1, *_2, *_3 \neq 0$. One can apply $M_{32}^a$ and $M_{43}^b$ for
appropriate choices of $a$, $b \in k$ to get
$\begin{bmatrix} 0 & 0 & * & * & * & *\end{bmatrix}$.
\item One of $*_1, *_2, *_3 = 0$. Use $T_{23}$ and $T_{34}$ to
get $*_1 = 0$. Now apply $M_{43}^a$ to get $\begin{bmatrix} 0 & 0 & * & * & * & *\end{bmatrix}$.
\item Two or more of $*_1, *_2, *_3 = 0$. Use $T_{23}$ and
$T_{34}$ to get $\begin{bmatrix} 0 & 0 & * & * & * & *\end{bmatrix}$.
\end{itemize}
We now have a basis vector of the form $\begin{bmatrix} 0 & 0 & * & * & * & *\end{bmatrix}$. We focus attention on
$*_3$ and $*_4$. We have the following cases.
\begin{itemize}
\item $*_3 = 0$: All spaces generated by vectors of the form
$\begin{bmatrix} 0 & 0 & 0 & * & * & *\end{bmatrix}$ are equivalent to
$\begin{bmatrix} 0 & 0 & 0 & 0 & 0 & 1\end{bmatrix}$ by applying
$M_{3 2}^a$, $M_{2 3}^b$, $M_{4 3}^c$ and $M_{3 4}^d$ for appropriate
$a$, $b$, $c$,  $d \in k^\times$.
\item $*_4 = 0$: All spaces generated by vectors of the form
$\begin{bmatrix} 0 & 0 & * & 0 & * & *\end{bmatrix}$  are equivalent to
$\begin{bmatrix} 0 & 0 & 0 & 0 & 0 & 1\end{bmatrix}$ by applying
 $M_{3 2}^a$, $M_{2 3}^b$, $M_{2 1}^c$ and $M_{1 2}^d$ for appropriate
$a$, $b$ , $c$, $d  \in k^\times$.
\item $*_3, *_4 \neq 0$: In this case we rescale the basis vector so that $*_3 = 1$ and
we have  $\begin{bmatrix} 0 & 0 & 1 & * & * & *\end{bmatrix}$.
One can show easily that the collection of $p^2$ vectors of
the form $\begin{bmatrix} 0 & 0 & 1 & y & * & *\end{bmatrix}$ with
$y \in k^\times$ form an orbit under the action of the
subgroup $\langle M_{1 2}^a, M_{2 3}^b \:|\: a, b \in k^\times\rangle$. One can thus restrict  to the case
$\begin{bmatrix} 0 & 0 & 1 & y & 0 & 0\end{bmatrix}$ with $y  \in k^\times$. The following chain of equivalences
\[
\begin{bmatrix} 0 & 0 & 1 & 1 & 0 & 0\end{bmatrix} \sim
\begin{bmatrix} 0 & 0 & 1 & 1 & y & 0\end{bmatrix} \sim
\begin{bmatrix} 0 & 0 & -y & 1 & y & 0\end{bmatrix} \sim \]
\[
\begin{bmatrix} 0 & 1 & 0 & 1 & y & 0\end{bmatrix} \sim
\begin{bmatrix} 0 & 0 & 1 & y & 1 & 0\end{bmatrix} \sim
\begin{bmatrix} 0 & 0 & 1 & y & 0 & 0\end{bmatrix}
\]
obtained by applying $M_{12}^y S_1^{-y} M_{2 1} T_{3 4} M_{1 2}^{-1}$ shows that
all vectors in this last case are equivalent to
$\begin{bmatrix} 0 & 0 & 1 & 1 & 0 & 0\end{bmatrix}$.
\end{itemize}
This completes the proof. Note that we  have not shown that
$\begin{bmatrix} 0 & 0 & 0 & 0 & 0 & 1\end{bmatrix}$  and
$\begin{bmatrix} 0 & 0 & 1 & 1 & 0 & 0\end{bmatrix}$  are
not equivalent to each other although this is in fact the case.
It can be deduced from our later results.
\end{proof}

We now use Theorem 1 to understand the action of $G$ on 2-dimensional spaces. When
specifying such a space we will give a pair of basis vectors in the form of a
$2 \times 6$ matrix. We have the following result.

\begin{theorem} Let $k^\times = \langle g \rangle$.
The space of 2-dimensional subspaces of $k^6$ form at most 4 orbits under the action of
$G$. In particular every 2-dimensional space is equivalent to one of the following.
\begin{enumerate}
\item \[
\begin{bmatrix}
0 & 0 & 0 & 0 & 0 & 1\\
1 & 0 & 0 & 0 & 0 & 0
\end{bmatrix}
\]
\item \[
\begin{bmatrix}
0 & 0 & 0 & 0 & 0 & 1\\
0 & 0 & 1 & 0 & 0 & 0
\end{bmatrix}
\]
\item \[
\begin{bmatrix}
0 & 0 & 0 & 0 & 0 & 1\\
0 & 0 & 1 & 1 & 0 & 0
\end{bmatrix}
\]
\item \[
\begin{bmatrix}
0 & 0 & 1 & 1 & 0 & 0\\
0 & 1 & 0 & 0 & g & 0
\end{bmatrix}
\]
\end{enumerate}
\end{theorem}
\begin{proof}
We select any 2-dimensional subspace U of $k^6$. Such a space contains $q + 1$
1-dimensional subspaces. By Theorem 1 each of these must be equivalent to
either
$\begin{bmatrix} 0 & 0 & 0 & 0 & 0 & 1\end{bmatrix}$  or
$\begin{bmatrix} 0 & 0 & 1 & 1 & 0 & 0\end{bmatrix}$.
We consider two cases. The first is where
U contains at least one subspace equivalent to
$\begin{bmatrix} 0 & 0 & 0 & 0 & 0 & 1\end{bmatrix}$
under the action of $G$ and the second is where it doesn't.

\medskip

\noindent
{\bf Case 1: U contains a subspace equivalent to
$\begin{bmatrix}  0 & 0 & 0 & 0 & 0 & 1 \end{bmatrix}$.}

\noindent
After applying a suitable element of $G$ we have
 \[
\begin{bmatrix}
0 & 0 & 0 & 0 & 0 & 1\\
* & * & * & * & * & *
\end{bmatrix}
\]
There are two subcases based on whether or not $*_1 = 0$. Before
we discuss these we list some elements of $G$ that stabilize the
first basis vector $\begin{bmatrix} 0 & 0 & 0 & 0 & 0 &
1\end{bmatrix}$. There are many such elements among the previously
listed generators of the group $G$ however there are four slightly
less obvious ones that will be useful. These include
$$
M_{32}^{b^{-1}} M_{23}^{-b} T_{23} S_2^{-b^{-1}} S_3^b,\quad M_{4
2}^{b^{-1}} M_{2 4}^{-b} T_{2 4}S_2^{-b^{-1}} S_4^b,
$$
$$
M_{31}^{b^{-1}} M_{1 3}^{-b} T_{1 3}S_1^{-b^{-1}} S_3^b,\quad M_{4
1}^{b^{-1}} M_{1 4}^{-b} T_{1 4}S_1^{-b^{-1}} S_4^b,
$$
where $b \in k^\times$. On a basis vector $R \:=\:
\begin{bmatrix} 1 & a_2 & a_3 & a_4 & a_5 & a_6\end{bmatrix}$ they
have the following effect.
\begin{description}
\item[(i)] $R M_{32}^{b^{-1}} M_{23}^{-b} T_{23} S_2^{-b^{-1}} S_3^b
= \begin{bmatrix} 1 & a_2 + b & a_3 & a_4 & a_5 & a_6 + b a_5\end{bmatrix}$
\item[(ii)] $R M_{4 2}^{b^{-1}} M_{2 4}^{-b} T_{2 4}S_2^{-b^{-1}} S_4^b
= \begin{bmatrix} 1 & a_2 & a_3 + b & a_4 & a_5 & a_6 - b a_4\end{bmatrix}$
\item[(iii)] $R M_{3 1}^{b^{-1}} M_{1 3}^{-b} T_{1 3}S_1^{-b^{-1}} S_3^b
= \begin{bmatrix} 1 & a_2 & a_3 & a_4 - b & a_5 & a_6 + b a_3\end{bmatrix}$
\item[(iv)] $R M_{4 1}^{b^{-1}} M_{1 4}^{-b} T_{1 4}S_1^{-b^{-1}} S_4^b
= \begin{bmatrix} 1 & a_2 & a_3 & a_4 & a_5 - b& a_6 - b a_2\end{bmatrix}$
\end{description}

We are now ready to consider the first subcase in which $*_1 \neq 0$.
First rescale so that $*_1 = 1$. Now apply the group elements described above that stabilize
the first basis vector to the second basis vector.  It
should be clear from the formulae above
that we can reduce all of the middle components to $0$. Subtracting
a multiple of the first vector to clear out the last component we see that the subspace
$U$ is equivalent to
 \[
\begin{bmatrix}
0 & 0 & 0 & 0 & 0 & 1\\
1 & 0 & 0 & 0 & 0 & 0
\end{bmatrix}
\]
This completes the first subcase.

In the second subcase we suppose that $*_1 = 0$. So that $U$ is given by
 \[
\begin{bmatrix}
0 & 0 & 0 & 0 & 0 & 1\\
0 & * & * & * & * & *
\end{bmatrix}
\]
As before we will make use of certain stabilizers of the first
basis vector. These are listed below together with their effects
on a vector of the form
$$
R = \begin{bmatrix} 0 & b_1 & b_2 & b_3 & b_4 & t \end{bmatrix}.
$$
\begin{description}
\item[(i)] $R T_{12}  =
\begin{bmatrix} 0 & b_3 & b_4 & b_1 & b_2 & t \end{bmatrix}$
\item[(ii)] $R T_{34}  =
\begin{bmatrix} 0 & b_2 & b_1 & b_4 & b_3 & t \end{bmatrix}$
\item[(iii)] $R M_{12}^a  =
\begin{bmatrix} 0 & b_1 & b_2 & b_3 + a b_1 & b_4 + a b_2 & t \end{bmatrix}$
\item[(iv)] $R M_{43}^a  =
\begin{bmatrix} 0 & b_1 + a b_2 & b_2 & b_3 + a b_4 & b_4 & t \end{bmatrix}$
\end{description}
\medskip

In the second basis vector at least one of the entries $*_2$,
$*_3$, $*_4$, $*_5$ must be nonzero. Using $T_{12}$ and $T_{34}$
we move this nonzero entry into the third position and then
rescale so that $*_3 = 1$ so that the second basis vector is now
of the form $\begin{bmatrix} 0 & * & 1 & * & * & *\end{bmatrix}$.
Applying $M_{4 3}^a$ and $M_{1 2}^b$ for appropriate choices of
$a$, $b \in k^\times$ we can get $*_2 = *_5 = 0$. We then subtract
a multiple of the first basis vector to get $*_6 = 0$ giving
 \[
\begin{bmatrix}
0 & 0 & 0 & 0 & 0 & 1\\
0 & 0 & 1 & * & 0 & 0
\end{bmatrix}
\]
Either $*_4 = 0$ or if it is nonzero it can be re-scaled so that
$*_4 = 1$ using $S_2^a$.

At the conclusion of Case 1 we see that any 2-dimensional subspace containing a
1-dimensional subspace equivalent to
$\begin{bmatrix} 0 & 0 & 0 & 0 & 0 & 1\end{bmatrix}$ must be equivalent to
one of the subspaces $(1)$, $(2)$ or $(3)$ listed in the statement of the Theorem.

\medskip

\noindent
{\bf Case 2: U does not contain a subspace equivalent to
$\begin{bmatrix} 0 & 0 & 0 & 0 & 0 & 1\end{bmatrix}$.}

\noindent
In this case all of the 1-dimensional subspaces must be equivalent to
$$
\begin{bmatrix} 0 & 0 & 1 & 1 & 0 & 0\end{bmatrix}.
$$
Fixing a basis we can apply an element of $G$ to reach a subspace of the form
\[
\begin{bmatrix}
0 & 0 & 1 & 1 & 0 & 0\\
* & * & * & * & * & *
\end{bmatrix}
\]
If $*_1 = *_2 = 0$ then subtracting a multiple of the first basis vector we see that
we can assume $*_3 = 0$. However a nonzero vector of the form
$\begin{bmatrix} 0 & 0 & 0 & * & * & *\end{bmatrix}$
is equivalent to
$\begin{bmatrix} 0 & 0 & 0 & 0 & 0 & 1\end{bmatrix}$
(see the proof of Theorem 1) and we have already
considered spaces $U$ containing such subspaces in Case 1. We can therefore assume that at least one of $*_1$ or $*_2$ is nonzero.

If $*_2 = 0$ then we switch $*_2$ and $*_1$ with $T_{23} S_2^{-1}$ (an element that
leaves the first basis vector unchanged). We rescale so that $*_2 = 1$ and
apply $M_{32}^a$ and $M_{1 4}^b$ for appropriate $a$, $b \in k^\times$ to
get $*_1 = 0$ and $*_6 = 0$. Subtracting a multiple of
the first basis vector to get $*_3 = 0$ we have now reduced to the case of a subspace of the form
\[
\begin{bmatrix}
0 & 0 & 1 & 1 & 0 & 0\\
0 & 1 & 0 & * & * & 0
\end{bmatrix}
\]
At least one of $*_4$ or $*_5$ must be nonzero otherwise the second basis vector
is equivalent to
$\begin{bmatrix} 0 & 0 & 0 & 0 & 0 & 1\end{bmatrix}$ and we are back in Case 1.
We make a note of the
following transformation.
\[
\begin{bmatrix} 0 & x & 0 & y & z & 0\end{bmatrix} S_2^\alpha S_4^\alpha \:=\:
\begin{bmatrix} 0 & x & 0 & \alpha y & \alpha^2 z & 0\end{bmatrix}
\qquad \qquad (\alpha \neq 0)  \]
This transformation leaves the $1$-dimensional subspace defined
by the first basis vector invariant. Applying this it is clear that every space with $*_4 = 0$
is equivalent to one of two spaces depending on whether or not $*_5$ is a square or
nonsquare element of $k$. Indeed we have
\[
\begin{bmatrix}
0 & 0 & 1 & 1 & 0 & 0\\
0 & 1 & 0 & 0 & * & 0
\end{bmatrix}
\:\sim\:
\begin{bmatrix}
0 & 0 & 1 & 1 & 0 & 0\\
0 & 1 & 0 & 0 & 1 & 0
\end{bmatrix}
\quad \text{or} \quad
\begin{bmatrix}
0 & 0 & 1 & 1 & 0 & 0\\
0 & 1 & 0 & 0 & g & 0
\end{bmatrix}
\]
where $k^\times = \langle g \rangle$.

We are thus left to consider the case where the second vector has
the form
$\begin{bmatrix} 0 & 1 & 0 & * & * & 0\end{bmatrix}$ and $*_4 \neq 0$.
We may assume that $*_4 = 1$ by applying $S_2^\alpha S_4^\alpha$ for suitable $\alpha$.

The following chain of equivalences is critical and relates the cases $*_4 = 0$ and $*_4 \neq 0$.
The transformations involved do not fix the first subspace so we record both.
Starting with
\[
\begin{bmatrix}
0 & 0 & 1 & 1 & 0 & 0\\
0 & 1 & 0 & 0 & t & 0
\end{bmatrix}
\]
add the first basis vector to the second to get
\[
\begin{bmatrix}
0 & 0 & 1 & 1 & 0 & 0\\
0 & 1 & 1 & 1 & t & 0
\end{bmatrix}
\]
and then apply $M_{43}^{-1} M_{34}^1 M_{12}^{-1}$
\[
\begin{bmatrix}
0 & -1 & 0 & 2 & 1 & 0\\
0 & 0 & 1 & 1-t & 0 & 0
\end{bmatrix}
\]
We now suppose that $t \neq 1$ and apply $S_1^\alpha$ with $\alpha = 1 - t$.
\[
\begin{bmatrix}
0 & t-1 & 0 & 2 & 1 & 0\\
0 & 0 & 1-t & 1-t & 0 & 0
\end{bmatrix}
\:\sim\:
\begin{bmatrix}
0 & 0 & 1 & 1 & 0 & 0 \\
0 & 1 & 0 & 2/(t -1) & 1/(t - 1) & 0
\end{bmatrix}
\]
The last equivalence results from re-scaling the vectors and
switching their order. To finish apply $S_2^\alpha S_4^\alpha$
with $\alpha = (t-1)/2$. After re-scaling the first basis vector
this gives
\[
\begin{bmatrix}
0 & 0 & 1 & 1 & 0 & 0 \\
0 & 1 & 0 & 1 & (t - 1)/4 & 0
\end{bmatrix}
\]

In summary if $x = (t-1)/4$ then
\[
\begin{bmatrix}
0 & 0 & 1 & 1 & 0 & 0 \\
0 & 1 & 0 & 1 & x & 0
\end{bmatrix}
\:\sim\:
\begin{bmatrix}
0 & 0 & 1 & 1 & 0 & 0\\
0 & 1 & 0 & 0 & t & 0
\end{bmatrix}
\]
provided that $t \neq 1$ or equivalently that $x \neq 0$. We are thus done provided
$x \neq 0$ since we have already considered matrices of the second form.
If $x = 0$ then the second basis vector
$\begin{bmatrix} 0 & 1 & 0 & 1 & 0 & 0\end{bmatrix}$ is equivalent to
$\begin{bmatrix} 0 & 0 & 0 & 1 & 0 & -1\end{bmatrix}$ under $T_{1 4}$ and this is equivalent to
$\begin{bmatrix} 0 & 0 & 0 & 0 & 0 & 1\end{bmatrix}$  putting
us back in Case 1.

So far we have shown in Case 2 that any 2-dimensional subspace (not already
covered by Case 1) is equivalent to
\[
\begin{bmatrix}
0 & 0 & 1 & 1 & 0 & 0 \\
0 & 1 & 0 & 0 & 1 & 0
\end{bmatrix}
\quad\text{or}\quad
\begin{bmatrix}
0 & 0 & 1 & 1 & 0 & 0\\
0 & 1 & 0 & 0 & g & 0
\end{bmatrix}
\]
where $k^\times = \langle g \rangle$. We now show that the first of these is also equivalent to a space
of the type covered in Case 1. Starting with
\[
\begin{bmatrix}
0 & 0 & 1 & 1 & 0 & 0 \\
0 & 1 & 0 & 0 & 1 & 0
\end{bmatrix}
\]
apply $M_{42}^1 M_{21}^{-1} M_{13}^1 T_{12} M_{24}^{-2} M_{42}^1 M_{21}^1 M_{32}^1$ to
get
\[
\begin{bmatrix}
1 & 0 & 3 & 0 & 0 & -1 \\
1 & 0 & 0 & 0 & 0 & 1
\end{bmatrix}
\]
Subtracting the first basis vector from the second we obtain
$\begin{bmatrix} 0 & 0 & -3 & 0 & 0 & 2\end{bmatrix}$ which is of the form
$\begin{bmatrix} 0 & 0 & * & 0 & * & *\end{bmatrix}$. All such nonzero vectors were shown to be equivalent to $\begin{bmatrix} 0 & 0 & 0 & 0 & 0 & 1\end{bmatrix}$
in the proof of Theorem 1. It follows that our 2-dimensional space is equivalent to
one of the three possibilities that arose in Case~1.
\end{proof}

\begin{definition}
We say that two quadratic presentations are {\em equivalent} if the associated subspaces
lie in the same orbit under the action of the group $G$.
\end{definition}

Theorem 2 implies that there are at most 4 types of quadratic
presentation up to equivalence. We now show that there are exactly
4. The main step is the following result which shows that the
presentations associated to $(1)$ and $(4)$ are not equivalent.

\begin{theorem}
Let $k^\times =  \langle g \rangle$. The Lie algebras $L_1 = L/\mathfrak{r}_1$ and $L_2 = L/\mathfrak{r}_2$
where
\[  \mathfrak{r}_1 =  \langle x_{13}, x_{14}, x_{23}, x_{24} \rangle \]
\[ \mathfrak{r}_2 =  \langle x_{12} , x_{34}, x_{14} - x_{23}, g x_{13} - x_{24} \rangle   \]
are not isomorphic.
\end{theorem}
\begin{proof}
If $L_1$ and $L_2$ were isomorphic then we would have induced isomorphisms on
their quotients by terms in the lower central series.  In particular we would have
$K_1 \cong K_2$ where $K_i = L_i/[[L_i,L_i],L_i]$ for $i = 1$, $2$.
The Lie algebra $K_1$ has several
elements with centralizer of dimension $5$. If
$K_1 \cong K_2$ this would imply that $K_2$ should also have $5$-dimensional element
stabilizers.
We will show this is not possible.

First write let us rewrite the condition $[v,w] = 0$ for $v$ and $w$ in $K_2$.
We start by writing
\begin{eqnarray*}
v &=&  a_1 x_1 + a_2 x_2 + a_3 x_3 + a_4 x_4 + a_{13} x_{13} + a_{14} x_{14}\\
w &=& b_1 x_1 + b_2 x_2 + b_3 x_3 + b_4 x_4 + b_{13} x_{13} + b_{14} x_{14}
\end{eqnarray*}
for some constants $a_i$, $a_{ij}$, $b_i$ and $b_{ij}$ in $k$. We have abused notation
slightly by using $x_i$ and $x_{ij}$ to also represent their  images in the quotient $K_2$.
The relations in
$L_2$ (and hence also in $K_2$) imply that  $x_{12} = x_{34} = 0$,
 $x_{23} = x_{14}$ and $x_{24} = g x_{13}$. We can now compute the bracket $[v,w] $
and simplify to get the equation
\[ 0 = [a_1b_3  + g^{-1}a_2 b_4  - a_3 b_1 - g^{-1} a_4 b_2] x_{13}
+ [a_1 b_4  - a_2 b_3 - a_3 b_2 - a_4 b_1]x_{14}. \]
We note that the constants $a_{ij}$ and $b_{ij}$ have completely disappeared
and play no further role in the argument.
Since $x_{13}$ and $x_{14}$ are linearly independent the coefficients must be zero so
we get a pair of equations.
Let $\mathbf{a}$ be the row vector $[a_1,a_2,a_3,a_4]$ and similarly for $\mathbf{b}$.
We will use the notation $M^T$ for the transpose of a vector or matrix. We have
two equations $\mathbf{a} M_1 \mathbf{b}^T = 0$ and $\mathbf{a} M_2 \mathbf{b}^T = 0$
where
\[
M_1 =
\begin{bmatrix}
0 & 0 & 1 & 0 \\
0 & 0 & 0 & g^{-1} \\
-1 & 0 & 0 & 0 \\
0 & -g^{-1} & 0 & 0
\end{bmatrix}
\qquad\text{and}\qquad
M_2 =
\begin{bmatrix}
0 & 0 & 0 & 1 \\
0 & 0 & 1 & 0 \\
0 & -1 & 0 & 0 \\
-1 & 0 & 0 & 0
\end{bmatrix}
.\]
If we let $\mathbf{c} = \mathbf{a} M_2$ then these equations become $\mathbf{c} \mathbf{b}^T = 0$
and $\mathbf{c} M \mathbf{b}^T = 0$ where $M = M_2^{-1} M_1$. For a given $\mathbf{c} \neq 0$
the space of solutions $\mathbf{b}$ has dimension $\leq 5$ with equality if and only if
$\mathbf{c} M = \lambda \mathbf{c}$ for some $\lambda \in k$. But this cannot happen since the matrix $M$
has eigenvalues $\pm \sqrt{g}$ which do not lie in $k$.
\end{proof}

\begin{corollary}
There are exactly $4$ orbits of $2$-dimensional (or $4$-dimensional) subspaces of $k^6$
under the action of $G$. Two of the associated presentations are mild and two are not mild.
There are exactly $2$ orbits of $1$-dimensional subspaces.
\end{corollary}
\begin{proof}
The orbits $(1)$ , $(2)$ and $(3)$ in Theorem 2 give rise to the mild (cycle)
presentation and two non-mild presentations. The corresponding
Lie algebra  $\mathfrak{g} = L/ \mathfrak{r} = \oplus_{n = 1}^\infty
\mathfrak{g}_n$ in each case can be distinguished from the others simply by computing
the dimension $a_n = \dim \mathfrak{g}_n$ for $n \leq 4$. Indeed we have
$a_1 = 4$, $a_ 2 = 2$ and
then
\begin{itemize}
\item[(1)] $a_3 = 4$, $a_4 = 6$;
\item[(2)] $a_3 = 5$;
\item[(3)] $a_3 = 4$, $a_4 = 7$.
\end{itemize}
It follows that the three associated presentations are not
equivalent and give rise to distinct orbits under the action of
$G$.

Theorem 3 shows that $(1)$ and $(4)$ are not equivalent, however observe that
they become equivalent over the extension
$k(\sqrt{g})$ and hence have the same Poincar\'e series over both $k(\sqrt{g})$ and
$k$. This is sufficient by \cite{La}, Proposition 3.2, to show that the presentation associated to $(4)$ is
strongly free and that $(4)$ cannot be equivalent to $(2)$ or $(3)$.

The statement about the $1$-dimensional subspaces now follows since if there
were only one orbit then the proof of Theorem 2 would yield
an upper bound of 3 on the number of $2$-dimensional subspaces.
\end{proof}

The arguments in the proof of Corollary 1 also yield the
following simple criterion for mildness.

\begin{corollary}
Let $G$ be a $4$-generated pro-$p$-group whose associated Lie
algebra $\mathfrak{g}=L/\mathfrak{r}$ over $\mathbb{F}_p$ is
quadratic of relation rank $4$. Then $G$ is mild if and only if
$\dim \mathfrak{g}_3  = 4$ and $\dim
\mathfrak{g}_4 = 6$.
\end{corollary}

One can find examples of Koch presentations which belong in each of the
four orbits described in Theorem 2. Indeed one can find $p$ and $S$
such that this is the case for the Galois group $G_S(p)$. We list
four such examples (the numberings are matched). In each case take $p = 3$.
\begin{itemize}
\item[(1)] $S = \{ 31, 37, 43, 67 \}$.
\item[(2)] $S = \{ 67, 79, 97, 127\}$.
\item[(3)] $S = \{ 61, 73, 79, 97 \}$.
\item[(4)] $S = \{ 31, 37, 61, 67 \}$.
\end{itemize}

\section{Open Questions}

Let $\rho_1,\dots,\rho_m$ be quadratic Lie polynomials in $m$
variables $x_1,\ldots,x_m$.
\begin{enumerate}[{\rm (A)}]
\item Find an algorithm for determining the strong freeness of
$\rho_1,\dots,\rho_m$.
\medskip

\item Find the number of inequivalent strongly free sequences $\rho_1,\dots,\rho_m$.
\medskip

\item If $\rho_1,\dots,\rho_m$ is strongly free, is it equivalent to one of the
form
$$
[x_1,x_2],[x_2,x_3],\ldots,[x_{m-1},x_m],[x_m,x_1]
$$
over an algebraically closed field?
\medskip

\item If $a_n$ is the dimension of the $n$-th homogeneous component
of $\mathfrak{g}=L/(\rho_1,\dots,\rho_m)$, is $\rho_1,\dots,\rho_m$
strongly free if
$$
\prod_{n\ge0}(1-t^n)^{a_n}=1-mt+mt^2 \text{ mod } t^c?
$$
for some $c$ depending only on $m$? Is $c=5$?
\medskip

\item In \cite{Sc}, Schmidt shows that under certain conditions (see \cite{Sc}, Theorem 2.1),
$\text{cd}(G_S(p))=2$ if $T\subset S$ and $G_T(p)$ is mild. Is
$G_S(p)$ mild under these conditions?
\end{enumerate}

\end{document}